# INTERACTING MULTI-CLASS TRANSMISSIONS IN LARGE STOCHASTIC NETWORKS

By Carl Graham and Philippe Robert

*UMR 7641 CNRS—École Polytechnique and INRIA Paris—Rocquencourt*

The mean-field limit of a Markovian model describing the interaction of several classes of permanent connections in a network is analyzed. Each of the connections has a self-adaptive behavior in that its transmission rate along its route depends on the level of congestion of the nodes of the route. Since several classes of connections going through the nodes of the network are considered, an original mean-field result in a multi-class context is established. It is shown that, as the number of connections goes to infinity, the behavior of the different classes of connections can be represented by the solution of an unusual nonlinear stochastic differential equation depending not only on the sample paths of the process, but also on its *distribution*. Existence and uniqueness results for the solutions of these equations are derived. Properties of their invariant distributions are investigated and it is shown that, under some natural assumptions, they are determined by the solutions of a fixed-point equation in a finite-dimensional space.

**1. Introduction.** Mathematical modeling of data transmission in communication networks has been the subject of an intense activity for some time now. For data transmission, the Internet network can be described as a very large distributed system with self-adaptive capabilities to the different congestion events that regularly occur at its numerous nodes. Various approaches have been used in this respect: control theory, ordinary differential equations, Markov processes, optimization techniques and so on. For data transmission one can distinguish two levels of mathematical abstraction:

(1) *Packet level.* For this approach, the self-adaptive behavior of the transmission protocols is analyzed. Packets are sent as long as no loss is detected and throughput grows linearly during that time. On the contrary

---

when a loss occurs, the throughput is sharply reduced by a multiplication by some factor $r < 1$ (the constant $1/2$ in general). This scheme is known as an Additive Increase and Multiplicative Decrease algorithm (AIMD). The congestion avoidance algorithm of the TCP protocol is an example of such a scheme.

Ott, Kemperman and Mathis [27] is one of the early mathematical studies. In this setting, after a scaling with respect to the loss rate, Dumas, Guillemin and Robert [11] proves various limit theorems for the resulting processes. A *fluid limit description* of the stochastic evolution of the throughput of a connection can then be obtained through a Markov process whose infinitesimal generator is given by

$$\Omega(f)(x) = af'(x) + bx(f(rx) - f(x))$$

for $f$ a $C^1$-function from $\mathbb{R}_+$ to $\mathbb{R}$ where $a$, $b > 0$ and $r \in (0, 1)$. The constant $a$ is the rate of linear increase and $r$ is the multiplicative factor. A precise formulation of this procedure is detailed in Section 2. See Dumas, Guillemin and Robert [11] for a survey of the literature in this domain. See also Chafai, Malrieu and Paroux [9].

(2) *Flow level.* It consists of looking at the state of the network on a quite long time scale so that the self-adaptive procedures of data transmission in the Internet can be described as a variation of a processor-sharing discipline at the level of the network. A connection is seen as a fluid characterized by the throughput it achieves given the other competing flows. A possible mathematical formulation is via an optimization problem: given $J$ classes of connections, when there are $x_j$ connections of class $j \in \{1, \dots, J\}$, the total throughput achieved by the connections of class $j$ is given by $\lambda_j$ so that the vector $(\lambda_j)$ is a solution of

$$\max_{\lambda \in \Lambda} \sum_{j=1}^{J} x_j U_j(\lambda_j / x_j),$$

where $\Lambda$ is the set of admissible throughputs which takes into account the capacity constraints of the network. The functions $(U_j)$ are defined as *utility* functions; various expressions have been proposed for them. See Kelly, Maulloo and Tan [21], Massoulié [23] and Massoulié and Roberts [24]. Kelly and Williams [22] gives a heavy-traffic analysis of a network in such a context. See also Srikant [31] for a more detailed overview of these models.

In this paper, the first approach, the fluid limit description of the stochastic evolution of the throughput of a connection, will be used to analyze the coexistence of numerous connections in a network with a general number of nodes. Up to now the stochastic models analyzed rigorously have considered a *single* node carrying *similar* flows. In this context, some approximated mean field analyses have been carried out for a single buffer receiving packets of a large number of connections (see, e.g., Adjih, Jacquet and Vvedenskaya [1] and Baccelli, McDonald and Reynier [3]). In particular, few rigorous, explicit, results have been obtained for models with several types of connections (e.g., relating different couples of nodes) and describing the evolution of the state of the whole network. In the general context of stochastic networks, mean-field approaches have been used in various situations notably by Dobrushin and his co-authors (see, e.g., Karpelevich, Pechersky and Suhov [20]).

The simultaneous transmission of several classes of permanent connections is analyzed here. A connection is described, via a fluid representation, introduced in Section 2 by its instantaneous throughput. A class of connections is characterized by the set of nodes it uses and by the rate of increase of its throughput when no loss occurs. The connections interact as follows: the loss rate of a given connection depends on the congestion of the nodes it traverses. If a loss occurs for a connection, its throughput is multiplied by some factor less than 1. In particular, the throughput of a connection depends on the state of *all* the nodes it is using and therefore on the other connections traversing any of these nodes.

*Mean-field analysis of multi-class systems.*    The evolution of a *large* number of connections divided into homogeneous classes will be considered. As it will be seen, a Markovian multi-class mean-field analysis has to be developed to investigate the properties of these networks. In this setting, chaoticity (i.e., convergence in law to i.i.d. random variables) is a classical notion. For statistically indistinguishable objects, modeled by systems of exchangeable random variables, chaoticity is equivalent to the fact that the empirical measures satisfy a weak law of large numbers. Sznitman [33] develops compactness–uniqueness methods, as well as coupling methods between exchangeable systems and an i.i.d. system which give stronger results. See also Méléard [26] and Graham [12, 13].

Many phenomena involve dissimilar objects, called particles, classified in a finite number of types, particles of a class being similar and *numerous* (see, e.g., Cercignani, Illner and Pulvirenti [8] and the review papers Bellomo and Stöcker [4], Graham [13], Grunfeld [16], Struckmaier [32]). Few convergence proofs for such multi-class systems exist, and most references directly use extensions of the limit equations for systems of identical particles without further justification. Some coupling techniques, such as those in Graham and Méléard [15] and Graham [13], using representations of the past history by random graphs and trees, may apply to some multi-class systems, but the models considered here do not have the specific structure required. Graham [14] studies natural notions of multi-exchangeability and chaoticity, and one of its scopes is to enable the adaptation of Sznitman's

compactness–uniqueness methods to multi-class systems, but these methods will not be used here.

*Nonlinear stochastic differential equations.* As it will be seen, the sample path of the evolution of the fluid representation of the throughput of a "typical" connection is converging to the solution of an unusual nonlinear SDE (stochastic differential equations) driven by Poisson processes. In a one-dimensional diffusion setting, the corresponding equation can be written as

$$dX(t) = \sigma[X(t), E(X(t))]\mathcal{N}(dt) + m[X(t), E(X(t))]\,dt, \qquad t \geq 0,$$

where $(\mathcal{N}(t))$ is a Poisson process on $\mathbb{R}_+$ and $\sigma$ and $m$ are sufficiently regular functions. If $\phi(t) = \mathbb{E}(X(t))$ were some fixed function, the above equation would be a standard (nonhomogeneous) SDE with a solution $(X(t))$ such that the value of $X(T)$ is a functional of the sample path of $(\mathcal{N}(s), s \leq T)$ and of $(\phi(s), s \leq T)$. But, since $\phi$ is dependent on the *distribution* of $X(t)$, the problem here is completely different. In a diffusion framework, existence results have been proved in Jourdain, Méléard and Woyczynski [19] and Benachour et al. [6] (see also Joudain and Malrieu [18] for further references and properties). In our case, existence and uniqueness results for the solution of the corresponding nonlinear SDE when Poisson processes are proved in Section 4. Moreover, in order to have sufficiently general initial conditions, that is, with a Gaussian tail distribution, for example, an important exponential estimate, relation (4.7), for these solutions is established. This inequality is a key ingredient to prove the main mean field convergence result and for the applications of the last section. Note that with deterministic initial conditions (or with a bounded support), the corresponding results would have been more straightforward to derive.

*Outline of the paper.* The mean-field limit of a Markovian model describing the interaction of $K$ classes of permanent connections in a network with $J$ nodes is analyzed. Each of the connections has a self-adaptive behavior in that its transmission rate along its route depends on the level of congestion of the nodes of the route. More formally, the state of a class $k$ connection is given by its throughput, that is, some nonnegative number $w_k \in \mathbb{R}_+$, and the quantity $u_j$ of a node $j \in \{1, \ldots, J\}$ is given as a weighted sum of all throughputs of the connections using this node. If $w$ is the state of a class $k$ connection, it increases linearly with slope $a_k(w, u)$ and jumps to $r_k w$ at rate $b_k(w, u)$ where $u = (u_j)$ is the weighted throughput vector and $r \in (0, 1)$ is the multiplicative factor. A sketch of the stochastic evolution of the $n$th class $k$ connection, $1 \leq n \leq N_k$, can be given by a fluid description through the ordinary differential equation

$$\frac{d}{dt} w_{k,n}(t) = a_k(w_{k,n}(t), u(t)) - (1 - r_k)w_{k,n}(t)b_k(w_{k,n}(t), u(t)),$$

where $u(t) = (u_j(t))$ and, for $1 \leq j \leq J$,

$$u_j(t) = \sum_{k=1}^{K} A_{jk} \sum_{p=1}^{N_k} w_{k,p}(t)$$

(see Section 3 for a more formal description). Related models are considered by Kelly, Maulloo and Tan [21].

Since several classes of connections going through the nodes of the network are considered, an original mean-field result in a multi-class context is established in the following context: for any $1 \leq k \leq K$, the number $N_k$ of class $k$ connections goes to infinity such that $N_k/(N_1 + \cdots + N_K)$ goes to $p_k$. Due to the lack of symmetry, an extension to multi-class systems of the notion of chaoticity and of Sznitman's coupling methods is developed. It is shown that, in the mean-field limit, the asymptotic behavior of a "typical" class $k$ connection, $1 \leq k \leq K$, can be represented by the solution of an unusual nonlinear stochastic differential equation $((W_k(t)), 1 \leq k \leq K)$ depending not only on the sample paths of the process, but also on its *distribution*. The proof of the mean field convergence relies on the convergence of the vector $(u(t))$, and, more precisely, that the relation

$$\lim_{N_k \to +\infty} \mathbb{E}\left( \left| \frac{1}{N_k} \sum_{p=1}^{N_k} w_{k,p}^N(t) - E(W_k(t)) \right| \right) = 0$$

holds for any $1 \leq k \leq K$ where the upper index $N$ refers to $(N_1, \ldots, N_K)$.

Existence and uniqueness results for the solutions $((W_k(t)), 1 \leq k \leq K)$ of these equations are established. Properties of their invariant distributions are investigated. It is shown, in particular, that for a quite large class of models the invariant distributions have an explicit expression given in terms of the solution of a fixed point equation in $\mathbb{R}_+^J$. Several examples of networks are discussed.

Section 2 recalls the main scaling results for the case of one class of connections in a one node network, and Section 3 defines the model of the networks analyzed in this paper. The properties of the limiting process associated with the behavior of a given transmission are investigated in Section 4. In Section 5, convenient notions of exchangeability and chaoticity in a multi-class system are introduced and mean-field results are obtained in this context. Section 6 is devoted to the analysis of possible invariant distributions of the limiting process.

## 2. A mathematical model of a single connection.
This section introduces the mathematical framework used to describe the evolution of the transmission rate of a single connection. In the following section, interactions of several connections within a network are presented and analyzed.

A connection between a source and a destination progressively increases its transmission rate until it receives some indication that the capacity along its path in the network is almost fully utilized. On the other hand, when congestion occurs, the throughput of the connection is drastically reduced. A given connection has a variable $W$ which gives the maximum number of packets that can be transmitted without receiving any acknowledgment from the destination. The variable $W$ is called the *congestion window size*. If all the $W$ packets are successfully transmitted, then $W$ is increased by 1 (additive part of the AIMD algorithm), so that $W + 1$ packets can be sent for the next round. Otherwise $W$ is divided by some factor (multiplicative part of the AIMD). An AIMD algorithm can thus be described as follows:

$$(2.1) \qquad W \to \begin{cases} W + 1, & \text{if no loss occurs among the } W \text{ packets,} \\ \lfloor rW \rfloor, & \text{otherwise,} \end{cases}$$

where $\lfloor x \rfloor$ is the integer part of $x \in \mathbb{R}$. The randomness for this model is given by the sequence of events when losses occur. One of the simplest and most reasonable models to consider is that the loss probability of a packet is given by $\varepsilon$ and that losses occur independently for each packet. This is the model analyzed in Ott, Kemperman and Mathis [27] at the fluid level and by Dumas, Guillemin and Robert [11]. The questionable independence property of losses is discussed thoroughly in Guillemin [17] where it is shown that this assumption is in fact pessimistic in that it underestimates the real efficiency of the algorithm.

*A fluid picture.* Under these assumptions, there is a natural Markov chain $(W_n^\varepsilon)$ whose transitions are given by (2.1). The variable $W_n^\varepsilon$ can be interpreted as the "throughput" at the beginning of the $(n+1)$th cycle. For a given $\varepsilon > 0$, little is known about the behavior of this Markov chain, about its invariant probability, for example. A scaling procedure has been developed in Dumas, Guillemin and Robert [11] consisting of looking at the continuous time Markov process

$$(\overline{W}_0^\varepsilon(t)) = (\sqrt{\varepsilon} W_{\lfloor t/\varepsilon \rfloor}^\varepsilon).$$

It is shown that, as $\varepsilon$ goes to 0, the family of Markov processes $(\overline{W}_0^\varepsilon(t), t \geq 0)$ converges in distribution to a Markov process $(\overline{W}_0(t), t \geq 0)$ whose infinitesimal generator $\Omega_0$ is given by

$$(2.2) \qquad \Omega_0(f)(x) = f'(x) + x(f(rx) - f(x)),$$

for $f$ a $C^1$-function from $\mathbb{R}_+$ to $\mathbb{R}$. The quantity $W_0(t)$ is interpreted as the instantaneous throughput of the connection at time $t$. Note that at this level, the "packet" description is not anymore valid since, in particular, one deal with a continuous state space. It is also proved that, as $\varepsilon$ goes to 0, the invariant distribution of $(\overline{W}_0^\varepsilon(t))$ converges to the invariant distribution of $(\overline{W}_0(t))$.

The linear growth rate $a$ of the throughput of the connection [the coefficient of $f'(x)$ which is 1 in (2.2)] is related to the inverse of RTT$^2$ where RTT is the round trip time, the duration of transmission between the two nodes (see Section 3). It should be noted that the RTT has another effect which is not taken into account in this paper; it delays the dynamic so that the evolution at time $t$ of a connection is a functional of the state of the network at time $s$, $s \in [t - \text{RTT}/2, t]$. See Mc Donald and Reynier [25] for an analysis of this phenomenon in the case of a network with one node.

In the following, a connection will be described as a Markov process similar to $(\overline{W}_0(t))$. Since different classes of connections will be considered, several parameters have to be introduced, and a connection will be represented as a Markov process $(W(t))$ whose infinitesimal generator is given by

$$(2.3) \qquad \Omega(f)(x) = af'(x) + bx(f(rx) - f(x))$$

for $f$ a $C^1$-function from $\mathbb{R}_+$ to $\mathbb{R}$. For $t \geq 0$, the quantity $W(t)$ should be thought as the instantaneous throughput of the connection at time $t$. As noted previously, its rate of increase $a$ is related to the round trip time between the source and the destination. The variable $b$ driving the loss rate will depend on the state of the nodes used by the connection. The quantities $a$, $b$ and $r$ will also depend on the class of the connection.

The Markov process $(W(t))$ increases linearly at rate $a$, and, given that $W(t) = x$, it jumps at time $t$ at rate $bx$ from $x$ to $rx$. As a stochastic process starting from $x \geq 0$, it can be represented as the solution $(W(t))$ of the following stochastic differential equation: $W(0) = x$ and

$$(2.4) \qquad \begin{aligned} dW(t) &= a\,dt - (1-r)W(t-)\mathcal{N}([0, bW(t-)], dt) \\ &= a\,dt - (1-r)W(t-)\int \mathbf{1}_{\{0 \leq x \leq bW(t-)\}}\mathcal{N}(dx, dt), \end{aligned}$$

where $\mathcal{N}$ is a Poisson process on $\mathbb{R}_+^2$ whose intensity measure is the Lebesgue measure on $\mathbb{R}_+^2$ and $W(t-)$ is the left limit of $W$ at $t$ (see Davis [10] and Robert [29]).

If $\rho = a/b$ and $(V(t)) = ((W(t)/\sqrt{\rho}))$, then the infinitesimal generator of $(V(t))$ is given by

$$\begin{aligned} \Omega_V(f)(x) &= \frac{a}{\sqrt{\rho}}f'(x) + b\sqrt{\rho}x(f(rx) - f(x)) \\ &= \sqrt{ab}(f'(x) + x(f(rx) - f(x))), \end{aligned}$$

which is proportional to the infinitesimal generator of (2.2). In particular, the invariant measure of $(V(t))$ and $(\overline{W}_0(t))$ are, therefore, the same. The density of the invariant distribution is given in the following proposition (see
Dumas, Guillemin and Robert [11] and Guillemin, Robert and Zwart [17]).

PROPOSITION 2.1. *If $W_0$ is a random variable with density*

$$(2.5) \quad H_{r,\rho}(x) = \frac{\sqrt{2\rho/\pi}}{\prod_{n=0}^{+\infty}(1 - r^{2n+1})} \sum_{n=0}^{+\infty} \frac{r^{-2n}}{\prod_{k=1}^{n}(1 - r^{-2k})} e^{-\rho r^{-2n}x^2/2}, \qquad x \geq 0,$$

*where $\rho = a/b$, then the solution of (2.4) with initial condition $W_0$ having
density $H_{r,\rho}$, is a stationary process and its expected value is given by*

$$(2.6) \qquad\qquad \mathbb{E}(W_0) = \sqrt{\rho}\,\psi(r)$$

*with*

$$\psi(r) = \sqrt{\frac{2}{\pi}} \prod_{n=1}^{+\infty} \frac{1 - r^{2n}}{1 - r^{2n-1}}.$$

At equilibrium, the quantity $\mathbb{E}(W_0/\sqrt{\varepsilon})$ can be seen as a first order ex-
pansion (in $\varepsilon$) of the asymptotic mean throughput of the connection.

**3. The stochastic model of the network.** The network has $J \geq 1$ nodes
and accommodates $K \geq 1$ classes of permanent connections. For $1 \leq k \leq K$,
the number of class $k$ connections is $N_k \geq 1$, one sets

$$N = (N_1, \ldots, N_K) \quad \text{and} \quad |N| = N_1 + \cdots + N_K.$$

An *allocation matrix* $A = (A_{jk}, 1 \leq j \leq J, 1 \leq k \leq K)$ with positive coef-
ficients describes the use of nodes by the connections. If $w_{n,k} \geq 0$ is the
throughput of the $n$th class $k$ connection, $1 \leq n \leq N_k$, the quantity $A_{jk}w_{n,k}$
is the weighted throughput at node $j$ of this connection. A simple example
would be to take $A_{jk} = 1$ or $0$, according to whether a class $k$ connection uses
node $j$ or not. The total weighted throughput $u_j$ of node $j$ by the various
connections is given by

$$u_j = \sum_{k=1}^{K} \sum_{n=1}^{N_k} A_{jk}w_{n,k}.$$

The quantity $u_j$ represents the level of congestion of node $j$, and, in par-
ticular, the loss rate of a connection going through it will depend on this
variable. A closely related model has been proposed in Kelly, Maulloo and
Tan [21] (see also Raghunathan and Kumar [28] in the context of wireless
network).

More precisely, for $1 \leq k \leq K$, the corresponding parameters $a$ and $b$ of
(2.3) for a class $k$ connection are given by functions $a_k : \mathbb{R}_+ \times \mathbb{R}_+^J \to \mathbb{R}_+$
and $b_k : \mathbb{R}_+ \times \mathbb{R}_+^J \to \mathbb{R}_+$, so that when the resource vector of the network is
$u = (u_j, 1 \leq j \leq J)$ and if the state of a class $k$ connection is $w_k$:

– its state increases linearly at rate $a_k(w_k, u)$;
– a loss for this connection occurs at rate $b_k(w_k, u)$, and in this case its state jumps to $r_k w_k$.

In view of (2.3), a natural form for these functions (with slight abuse of notation) is

$$(3.1) \qquad a_k(w_k, u) = a_k(u) \quad \text{and} \quad b_k(w_k, u) = w_k \beta_k(u).$$

For example, one can take

$$(3.2) \qquad \begin{cases} a_k(u) = \left( \tau_k + \displaystyle\sum_{j=1}^{J} t_{jk}(u_j) \right)^{-1}, \\ \beta_k(u) = \delta_k + \displaystyle\sum_{j=1}^{J} d_{jk}(u_j), \end{cases}$$

where $\tau_k > 0$ is the round trip time between source and destination and $\delta_k$ is the loss rate of class $k$ connections in a noncongested network. For $1 \leq j \leq J$, $t_{jk}(u_j)$ is their delay at node $j$ when the total weighted throughput is $u_j$. In particular, $t_{jk}(0) = 0$ and $t_{jk}(u) = 0$ when $A_{jk} = 0$. Similarly, $d_{jk}(u_j)$ is the loss rate at node $j$, and it is 0 when $A_{jk} = 0$.

*A Markovian representation.* The Markov process describing the state of the ongoing connections is given by, for $t \geq 0$,

$$W^N(t) = (W^N_{n,k}(t), 1 \leq n \leq N_k, 1 \leq k \leq K),$$

where $W^N_{n,k}(t)$ is the state of the $n$th flow of class $k$ at time $t$.

As for (2.4), it can be represented by the solution of the following stochastic differential equation: for $1 \leq k \leq K$ and $1 \leq n \leq N_k$,

$$\begin{aligned}(3.3) \quad dW^N_{n,k}(t) = {} & a_k(W^N_{n,k}(t-), U^N(t-)) \, dt \\ & - (1 - r_k) W^N_{n,k}(t-) \int \mathbf{1}_{\{0 \leq z \leq b_k(W^N_{n,k}(t-), U^N(t-))\}} \mathcal{N}_{n,k}(dz, dt)\end{aligned}$$

with $U^N(t) = (U^N_j(t), 1 \leq j \leq J)$ and

$$U^N_j(t) = \sum_{k=1}^{K} A_{jk} \sum_{i=1}^{N_k} W^N_{n,k}(t),$$

where $(\mathcal{N}_{n,k}, 1 \leq k \leq K, 1 \leq n \leq N_k)$ are independent Poisson processes on $\mathbb{R}^2_+$ whose intensity measure is a Lebesgue measure on $\mathbb{R}^2_+$. The following proposition is proved with standard arguments.

PROPOSITION 3.1. *If, for $1 \leq k \leq K$, the functions $a_k$ and $b_k$ are locally bounded, and $a_k$ is Lipschitz, then there is pathwise existence and uniqueness for a solution of stochastic differential equation (3.3).*

*Mean field scaling.* Due to its high-dimensional state space, the system of coupled stochastic differential equations (3.3) does not seem to be mathematically tractable as such. A scaling is used to investigate the qualitative properties of the corresponding Markov processes. It is assumed that for any $1 \leq k \leq K$, the number of class $k$ connections $N_k$ goes to infinity so that

$$(3.4) \qquad \frac{N_k}{|N|} = \frac{N_k}{N_1 + \cdots + N_K} \to p_k,$$

where $p_k \geq 0$ and $p_1 + \cdots + p_K = 1$.

The variable $U^N$ is accordingly scaled as $\overline{U}^N = U^N/|N|$ in the functions $a_k$ and $b_k$ so that (3.3) becomes the following: for $1 \leq k \leq K$ and $1 \leq n \leq N_k$,

$$(3.5) \qquad \begin{aligned} &dW_{n,k}^N(t) \\ &= a_k(W_{n,k}^N(t-), \overline{U}^N(t-)) \, dt \\ &\quad - (1 - r_k) W_{n,k}^N(t-) \int \mathbf{1}_{\{0 \leq z \leq b_k(W_{n,k}^N(t-), \overline{U}^N(t-))\}} \mathcal{N}_{n,k}(dz, dt) \end{aligned}$$

with $\overline{U}^N(t) = (\overline{U}_j^N(t), 1 \leq j \leq J)$ and

$$\overline{U}_j^N(t) \stackrel{\text{def}}{=} \frac{1}{N} U_j^N(t) = \sum_{k=1}^K \frac{N_k}{|N|} A_{jk} \overline{W}_k^N(t) \qquad \text{with } \overline{W}_k^N(t) = \frac{1}{N_k} \sum_{n=1}^{N_k} W_{n,k}^N(t).$$

Note that there is a difference in (3.5) with the initial equation (3.3) since a factor $1/N$ has been introduced in functions $a_k$ and $b_k$, $1 \leq k \leq K$, with the term $\overline{U}^N$. This is a scaling procedure in order to get a nontrivial scaling limit when $N$ gets large. This amounts to scale the capacity of routers with the parameter $N$. This is a multi-class mean-field system with an interaction through the *scaled weighted throughput* vector $(\overline{U}^N(t))$.

For $1 \leq k \leq K$, it is natural to introduce the *empirical measure* for the stochastic processes associated to class $k$ connections,

$$\Lambda_k^N = \frac{1}{N_k} \sum_{n=1}^{N_k} \delta_{(W_{n,k}^N(t), t \geq 0)},$$

where $\delta_{(x(t), t \geq 0)}$ is the Dirac mass at $(x(t), t \geq 0)$.

*Notation and conventions.* For $x = (x_m) \in \mathbb{R}^M$ for some $M \in \mathbb{N}$, $\|x\|$ denotes $\max(|x_m| : 1 \leq m \leq M)$ and if $z : \mathbb{R}_+ \to R^M$ and $T > 0$, one defines

$$\|z\|_T \stackrel{\text{def}}{=} \sup_{0 \leq s \leq T} \|z(s)\|.$$

Depending on the context the expressions

$$(z_m(t), 1 \leq m \leq M, t \geq 0) \quad \text{or} \quad (z_m(t), 1 \leq m \leq M)$$

are used for the function $(z(t))$.

If $H$ is a complete metric space, $\mathcal{D}(\mathbb{R}_+, H)$ denotes the Skorohod space of functions with values in $H$, continuous on the right and with left limits at any point of $\mathbb{R}_+$ (see Billingsley [7]). The variable $\Lambda_k^N$ has values in the set $\mathcal{P}(\mathcal{D}(\mathbb{R}_+, \mathbb{R}_+))$ of probability distributions on $\mathcal{D}(\mathbb{R}_+, \mathbb{R}_+)$.

## 4. Analysis of the nonlinear limiting process.

In view of (3.5), because of the symmetry properties within each class of connection, it is natural to expect a mean-field convergence phenomenon to hold when $N$ gets large. For $1 \leq k \leq K$, the evolution of a *typical* class $k$ connection should converge in law. More formally, the empirical measure $\Lambda_k^N$ *should* converge in law to the distribution of a stochastic process $(W(t)) = ([W_k(t), 1 \leq k \leq K], t \geq 0)$ solution of the nonlinear stochastic differential equation

$$
\begin{aligned}
dW_k(t) = {} & a_k(W_k(t-), u_W(t))\, dt \\
& - (1 - r_k) W_k(t-) \int \mathbf{1}_{\{0 \leq z \leq b_k(W_k(t-), u_W(t))\}} \mathcal{N}_k(dz, dt)
\end{aligned}
\tag{4.1}
$$

with $u_W(t) = (u_{W,j}(t), 1 \leq j \leq J)$ and, for $1 \leq j \leq J$,

$$
u_{W,j}(t) = \sum_{k=1}^K A_{jk} p_k \mathbb{E}(W_k(t)),
$$

where $(\mathcal{N}_k, 1 \leq k \leq K)$ are i.i.d. Poisson point processes on $\mathbb{R}_+^2$ with Lebesgue characteristic measure.

*A nonlinear stochastic differential equation.* In these equations, the interaction between coordinates of $(W(t))$ depends, in a nonlinear way, on the mean utilization vector $(u(t))$ which is a linear functional of mean values $(\mathbb{E}(W_k(t), 1 \leq k \leq K))$. In particular, the infinitesimal generator of the process $(W(t))$ depends *on the law* of $(W(t))$ and not only on its sample path up to time $t$ as it is usually the case.

The properties of the solutions to (4.1) are analyzed in this section. The convergence results of the mean-field scaling are the subject of Section 5.

*Nonlinear martingale problems.* Throughout this section, the natural filtration $(\mathcal{F}_t)$ used to investigate the solutions $(W(t))$ to (4.1) is defined as follows: for $t \geq 0$, $\mathcal{F}_t$ is the $\sigma$-field defined by

$$
\mathcal{F}_t = \sigma \langle \mathcal{N}_k(A \times B) : 1 \leq k \leq K, A \in \mathcal{B}(\mathbb{R}_+), B \in \mathcal{B}([0, t]) \rangle,
$$

where $\mathcal{B}(H)$ is the space of Borel sets of $H$ and $\mathcal{N}(A \times B)$ is the number of points of the point process $\mathcal{N}$ in $A \times B$. The elementary lemma will be used in the following.

Lemma. *If $(Y(t))$ and $(Z(t)) \in \mathcal{D}(\mathbb{R}_+, \mathbb{R}_+)$ are adapted processes to the filtration $(\mathcal{F}_t)$, then*

$$(I_Z(t)) = \left( \int_0^t \int Y(s-) \mathbf{1}_{\{0 \leq z \leq Z(s-)\}} [\mathcal{N}(dz, ds) - dz \, ds] \right)$$

$$= \left( \int_0^t Y(s-) [\mathcal{N}([0, Z(s-)], ds) - Z(s) \, ds] \right)$$

*is a local $(\mathcal{F}_t)$-martingale.*

Proof. By localization with a stopping time, $(Y(t))$ and $(Z(t))$ can be assumed to be uniformly bounded by some constant $C$.

If, for $T \geq 0$, $(Z(t), 0 \leq t \leq T)$ has a finite set $V$ of possible values, then, for $t \leq T$,

$$I_Z(t) = \sum_{v \in V} \int_0^t Y(s-) \mathbf{1}_{\{Z(s-) = v\}} (\mathcal{N}([0, v], ds) - v \, ds).$$

The classical properties of Poisson processes show that the martingale property holds on $[0, T]$ for each term of the sum in the above expression and consequently for the process $I_Z$.

Returning to the general case, since

$$\sup_{0 \leq t \leq T} |I_Z(t)| \leq C(\mathcal{N}([0, C] \times [0, T]) + CT)$$

by approximating $(Z(t), 0 \leq t \leq T)$ from below by processes with a finite number of possible values, and by using Lebesgue's theorem, one obtains the desired martingale property. □

If $W$ is a solution of (4.1), for $1 \leq k \leq K$, $t \geq 0$ and $f$ is a $C^1$-function from $\mathbb{R}_+$ to $\mathbb{R}$, by defining

$$
\begin{aligned}
(4.2) \quad M_k^f(t) = {}& f(W_k(t)) - f(W_k(0)) - \int_0^t a_k(W_k(s), u_W(s)) f'(W_k(s)) \, ds \\
& - \int_0^t b_k(W_k(s), u_W(s)) (f(r_k W_k(s)) - f(W_k(s))) \, ds,
\end{aligned}
$$

then $(M_k^f(t))$ is a local martingale. For $C^1$-functions $f$ and $g$ and $k \neq l$, the local martingales $M_k^f$ and $M_l^g$ are independent and the Doob–Meyer brackets are given by

$$\langle M_k^f \rangle_t = \int_0^t b_k(W_k(s), u_W(s)) (f(r_k W_k(s)) - f(W_k(s)))^2 \, ds$$

(see Rogers and Williams [30]).

A nonlinear martingale problem formulation of (4.1), is that there exists a probability $P$ on the space $\mathcal{D}(\mathbb{R}_+, \mathbb{R}_+^K)$ of sample paths $(W(t)) = ((W_k(t), 1 \leq k \leq K), t \geq 0)$, such that for any $C^1$-function $f$, the process $M_k^f$ defined by formula (4.2) is a local martingale for $P$.

*Existence and uniqueness results.* The following proposition is the central technical result used to establish the existence and uniqueness of solutions of (4.1), as well as the mean-field convergence result.

PROPOSITION 4.1.  *Let the functions $(a_k)$ be bounded, and $(a_k)$ and $(b_k)$ be Lipschitz. For any nonnegative continuous $\mathbb{R}_+^J$-valued function $u = (u(t))$ and any $\mathbb{R}_+^K$-valued random variable $X_0$, let $\phi(X_0, u) = (X_k(t))$ denote the solution of the*

$$
\begin{aligned}
(\mathcal{E}_u) \qquad dX_k(t) = {} & a_k(X_k(t-), u(t)) \, dt \\
& - (1 - r_k) X_k(t-) \\
& \times \int \mathbf{1}_{\{0 \leq z \leq b_k(X_k(t-), u(t))\}} \mathcal{N}_k(dz, dt), \qquad 1 \leq k \leq K,
\end{aligned}
$$

*starting at $(X_k(0)) = X_0$. Then, for any such functions $u$ and $u'$ and random variables $X_0$ and $X_0'$,*

$$
\begin{aligned}
(4.3) \qquad & \mathbb{E}[\|\phi(X_0, u) - \phi(X_0', u')\|_t \mid X(0), X'(0)] \\
& \leq [\|X_0 - X_0'\| + \int_0^t C(s) \|u - u'\|_s \, ds] e^{tC(t)}, \qquad t \geq 0,
\end{aligned}
$$

*where $C(t) = A(1 + t + \|u\|_t + \|X_0\| + \|X_0'\|)$ and $A$ is a constant depending only on the functions $(a_k)$ and $(b_k)$.*

PROOF.  The quantity $L_b$ [resp. $L_a$] denotes the maximum of the Lipschitz constants for the function $b_k$, [resp. $a_k$], $1 \leq k \leq K$. Denote by $(X_k(t)) = (\phi(X_0, u)(t))$ and $(X_k'(t)) = (\phi(X_0', u')(t))$.

Note that, since the functions $(a_k)$ are bounded, by some constant $\|a\|$, say, any solution $(Z(t))$ of $(\mathcal{E}_u)$ satisfies the relation

$$
(4.4) \qquad\qquad Z_k(t) \leq Z_k(0) + \|a\| t.
$$

The quantity $\mathbb{E}^0(\cdot)$ will refer to the conditional expectation given $X(0)$ and $X'(0)$. For $t \geq 0$,

$$
\begin{aligned}
(4.5) \qquad \|X_k - X_k'\|_t \\
\leq {} & |X_k(0) - X_k'(0)| \\
& + \int_0^t |a_k(X_k(s), u(s)) - a_k(X_k'(s), u'(s))| \, ds
\end{aligned}
$$

$$+ \int_0^t \int |X_k(s-)\mathbf{1}_{\{z \leq b_k(X_k(s-),u(s))\}}$$
$$- X_k'(s-)\mathbf{1}_{\{z \leq b_k(X_k'(s-),u'(s))\}}|\mathcal{N}(dz, ds).$$

By using the above lemma and the upper bound on $(X_k(t))$ which shows that the local martingale of this lemma is indeed a martingale, the conditional expected value of the last term of this expression is given by

$$\mathbb{E}^0\Big(\int_0^t \int |X_k(s)\mathbf{1}_{\{z \leq b_k(X_k(s),u(s))\}} - X_k'(s)\mathbf{1}_{\{z \leq b_k(X_k'(s),u'(s))\}}|\,dz\,ds\Big)$$

$$(4.6) \qquad \leq \mathbb{E}^0\Big(\int_0^t |X_k(s) - X_k'(s)|b_k(X_k(s), u(s))\,ds\Big)$$

$$+ \mathbb{E}^0\Big(\int_0^t X_k'(s)|b_k(X_k(s), u(s)) - b_k(X_k'(s), u'(s))|\,ds\Big).$$

By using the growth conditions and the Lipschitz properties of $(a_k)$ and $(b_k)$, one gets that

$$\mathbb{E}^0(\|X_k - X_k'\|_t) \leq |X_k(0) - X_k'(0)|$$
$$+ \int_0^t C_0(s)(\|u(s) - u'(s)\| + \mathbb{E}^0(\|X_k - X_k'\|_s))\,ds$$

with

$$C_0(t) = L_a + L_b(2\|a\|_t + \|u\|_t) + \sum_{k=1}^K b_k(0,0)$$

$$+ L_b[X_k(0) + X_k'(0)] + \sum_{j=1}^J L_b u_j(0)$$

so that $(C_0(t))$ can be replaced by $(C(t))$ as in the statement of the proposition. Gronwall's inequality gives the estimation

$$\mathbb{E}^0(\|X_k - X_k'\|_t) \leq \Big[|X_k(0) - X_k'(0)| + \int_0^t C(s)\|u(s) - u'(s)\|\,ds\Big]e^{\int_0^t C(s)\,ds}$$

$$\leq \Big[|X_k(0) - X_k'(0)| + \int_0^t C(s)\|u - u'\|_s\,ds\Big]e^{tC(t)}.$$

The proposition is complete.  $\square$

One of the main motivations of this study is to obtain results which are valid for functions $(a_k)$ and $(b_k)$ of the form (3.1) and (3.2) which have basically a quadratic behavior. To control the evolution of the vector $(W(t))$ and describe its stationary behavior, initial conditions cannot be assumed to be uniformly bounded; for this reason exponential and Gaussian moment assumptions are introduced.

CONDITION (C). It is said to hold for a family of random variables $\{X_0^\alpha, \alpha \in \mathcal{S}\}$ in $\mathbb{R}_+^K$, for $(b_k)$, and for $\varepsilon > 0$ when at least one of the two conditions is satisfied:

(1) For any $1 \leq k \leq K$, the function $b_k : \mathbb{R}_+ \times \mathbb{R}_+^J \to \mathbb{R}_+$ is Lipschitz and the variables $X_0^\alpha$, $\alpha \in \mathcal{S}$, have a uniform exponential moment of order $\varepsilon$,

$$\sup_{\alpha \in \mathcal{S}} \mathbb{E}(\exp(\varepsilon \|X_0^\alpha\|)) < \infty.$$

(2) For any $1 \leq k \leq K$, $b_k(w, u) = w \beta_k(u)$ and $\beta_k : \mathbb{R}_+^J \to \mathbb{R}_+$ is Lipschitz and the variables $X_0^\alpha$, $\alpha \in \mathcal{S}$, have a uniform Gaussian moment of order $\varepsilon$,

$$\sup_{\alpha \in \mathcal{S}} \mathbb{E}(\exp(\varepsilon \|X_0^\alpha\|^2)) < \infty.$$

The following theorem establishes the existence and uniqueness result of a solution of (4.1).

THEOREM 4.2. *If the functions* $a_k : \mathbb{R}_+ \times \mathbb{R}_+^J \to \mathbb{R}_+$, $1 \leq k \leq K$, *are bounded and Lipschitz, and if Condition* (C) *holds for* $W_0$, $(b_k)$ *and* $\varepsilon > 0$, *then there is pathwise existence and uniqueness of a solution* $(W(t))$ *of the nonlinear stochastic differential equation* (4.1) *starting at* $W_0$.

*In this case, the solution depends continuously on the initial condition in the following way: if* $(W(t))$ *and* $(W'(t))$ *are solutions of* (4.1) *with respective initial conditions* $W_0$ *and* $W_0'$, *having the same moment in Condition* (C), *for* $T \geq 0$, *there exists a constant* $A_T$ *such that*

$$(4.7) \quad \begin{aligned} &\mathbb{E}(\|W - W'\|_T) \\ &= \mathbb{E}\Big(\sup_{s \leq T} \|W(s) - W'(s)\|\Big) \\ &\leq A_T \mathbb{E}(\|W_0 - W_0'\| e^{\varepsilon(\|W_0\|^\ell + \|W_0'\|^\ell)/2}) e^{\mathbb{E}(\exp[\varepsilon(\|W_0\|^\ell + \|W_0'\|^\ell)])}, \end{aligned}$$

*where* $\ell = 1$ *or* $2$, *depending on the moment in Condition* (C), *whether it is exponential or Gaussian.*

PROOF. It is assumed for the moment that (1) of Condition (C) is satisfied. The functions $(a_k)$ are bounded by the constant $\|a\|$. Recall that for an integrable process $(Z(t)) = (Z_k(t))$, for $t \geq 0$,

$$u_Z(t) = (u_{Z,j}(t), 1 \leq j \leq J) = \left( \sum_{k=1}^K A_{jk} p_k \mathbb{E}(Z_k(t)), 1 \leq j \leq J \right).$$

Let $\alpha$ denote the constant $\max(A_{j1} p_1 + A_{j2} p_2 + \cdots + A_{jK} p_K : 1 \leq j \leq J)$.

The sequence $(W^n)$ is constructed by induction as follows:

– $W^0 \equiv W(0)$;

– for $n \geq 1$, $W^n$ is the solution of $(\mathcal{E}_u)$ (see Proposition 4.1) with $W^n(0) = W(0)$ and $(u(t)) = (u_{W^{n-1}}(t))$.

Note that by the growth condition $W^n(t) \leq W^n(0) + \|a\| t$ for $t \geq 0$ and $n \geq 1$, and consequently, $\|u_{W^n}\|_t \leq \alpha(\|W(0)\| + \|a\| t)$. With the notation of Proposition 4.1, $W^n = \phi(W(0), u_{W^{n-1}})$ for $n \geq 1$, so that (4.3) gives the relation

$$\mathbb{E}[\|W^{n+1} - W^n\|_t | W(0)] \leq \left[ \int_0^t C(s) \|u_{W^n} - u_{W^{n-1}}\|_s \, ds \right] e^{tC(t)}, \qquad t \geq 0,$$

with $C(t) = A_0(1 + t + \|W(0)\|)$ where $A_0$ is a constant depending only on $(a_k)$ and $(b_k)$. Because of the assumption on the existence of an exponential moment for $W(0)$, by taking $t_0 = \varepsilon/(2A_0)$, then $\delta = \alpha \mathbb{E}(C(t_0) \exp(t_0 C(t_0))) < +\infty$. One gets that for $t \leq t_0$,

$$\|u_{W^{n+1}} - u_{W^n}\|_t \leq \delta \int_0^t \|u_{W^n} - u_{W^{n-1}}\|_s \, ds, \qquad n \geq 1,$$

so that

$$\|u_{W^{n+1}} - u_{W^n}\|_t \leq \frac{(\delta t)^n}{n!} \|u_{W^1} - u_{W^0}\|_{t_0}, \qquad n \geq 1,$$

with relation (4.3). One gets, therefore, that

$$\mathbb{E}(\|W^{n+1} - W^n\|_{t_0}) \leq \frac{(\delta t_0)^n}{n!} \|u_{W^1} - u_{W^0}\|_{t_0}$$

holds almost surely. Consequently, the sequence of processes $(W^n(t), t \leq t_0)$ converges uniformly to some process $(W(t), t \leq t_0)$ and

$$\mathbb{E}(\|W^n - W\|_{t_0}) \leq \frac{(\delta t_0)^n}{n!} e^{\delta t_0} \|u_{W^1} - u_{W^0}\|_{t_0}.$$

The above relation also shows that on the interval $[0, t_0]$, the sequence of functions $(\mathbb{E}(W^n(t)))$ converges uniformly to $(\mathbb{E}(W(t)))$, and, in particular, $(u_{W^n}(t))$ converges uniformly to $(u_W(t))$ on this interval.

By definition of $W^{n+1}$, for $t \leq t_0$ and $1 \leq k \leq K$,

$$W_k^{n+1}(t) = W_k(0) + \int_0^t a_k(W_k^n(s), u_{W^n}(s)) \, ds$$

$$- (1 - r_k) W_k^n(t-) \int_0^t \mathcal{N}([0, b_k(W_k^n(s-), u_{W^n}(s))], ds),$$

the uniform convergences, the Lipschitz assumption for $a_k$ and $b_k$ and continuity properties of the Poisson process $\mathcal{N}_k$ show that, almost surely, the limit $(W(t))$ indeed satisfies relation (4.1) on $[0, t_0]$.

If $(W'(t))$ is another solution of (4.1) with $W'(0) = W(0)$, by using again relation (4.3), one gets that, for $t \leq t_0$,

$$\|u_W - u_{W'}\|_t \leq \delta \int_0^t \|u_W - u_{W'}\|_s \, ds$$

so that $u_W = u_{W'}$ on $[0, t_0]$, and, consequently, $W = W'$ on this interval.

Since $\|W(t_0)\| \leq \|W_0\| + \|a\| t_0$, Condition (C) holds for $W(t_0)$ and the same $\varepsilon > 0$, by taking it as an initial condition and by translating the Poisson processes by $t_0$ (on the second coordinate), one can show existence and uniqueness results on the time interval $[t_0, 2t_0]$, and therefore by induction on $\mathbb{R}_+$.

If $(W(t))$ and $(W'(t))$ are the solutions of (4.1) with respective initial conditions $W_0$ and $W_0'$, relation (4.3) gives that

$$\|u_W - u_{W'}\|_t \leq \alpha \mathbb{E}\!\left(\left[\|W_0 - W_0'\| + \int_0^t \|u_W - u_{W'}\|_s C(s) \, ds\right] e^{tC(t)}\right),$$

then, for $t < t_0$, by using again Gronwall's lemma, and the expression of $C(t)$, one gets that there exists a constant $A$ depending only on the functions $(a_k)$ and $(b_k)$ such that

$$\mathbb{E}(\|u_W - u_{W'}\|_t) \leq A \mathbb{E}(\|W_0 - W_0'\|) e^{\varepsilon(\|W_0\| + \|W_0'\|/2)} e^{\mathbb{E}(\exp(\varepsilon(\|W_0\| + \|W_0'\|)))}$$

for $t \leq t_0$. By plugging this inequality in relation (4.3), one gets a similar relation for $\mathbb{E}(\|W - W'\|_{t_0})$. With the same procedure, an analogous inequality can be obtained on the interval $[0, nt_0]$ for any $n \geq$ (with a constant $A$ depending on $n$). The theorem has therefore been proved on $\mathbb{R}_+$ in case (1) of Condition (C).

For the case (2) of Condition (C), when $b_k(w, u) = w\beta_k(u)$ and $\beta_k$ is Lipschitz, the corresponding inequality (4.6) in this setting is the relation

$$|X_k(s)\beta_k(u(s)) - X_k'(s)\beta_k(u'(s))|$$
$$\leq X_k(s)|\beta_k(u(s)) - \beta_k(u'(s))| + \beta_k(u'(s))|X_k(s) - X_k'(s)|$$

and the rest of the proof is analogous. The supplementary multiplication by $W_k(s)$ of the Lipschitz bounds is handled by the Gaussian moment assumption. □

## 5. The mean-field limit for converging initial data.

In this section, the convergence results to the limiting stochastic process analyzed in Section 4 are established. In the context of mean-field convergence, the notions of exchangeability and chaoticity play a fundamental role. If they are classical in the context of single-class systems (see, e.g., Aldous [2]), these properties have to be extended to investigate the stochastic model of the multi-class network considered in this paper.

DEFINITION.   The sequence of random variables $(X_{n,k}, 1 \leq n \leq N_k, 1 \leq k \leq K)$ is said to be *multi-exchangeable* if its law is invariant under permutation of the indexes *within* the classes: for $1 \leq k \leq K$ and any permutations $\sigma_k$ of $\{1, \ldots, N_k\}$, the equality in distribution

$$(X_{\sigma_k(n),k}, 1 \leq n \leq N_k, 1 \leq k \leq K) \overset{\text{dist}}{=} (X_{n,k}, 1 \leq n \leq N_k, 1 \leq k \leq K)$$

holds.

A sequence $(X_{n,k}^N, 1 \leq n \leq N_k, 1 \leq k \leq K)$ of multi-class random variables indexed by $N = (N_k) \in \mathbb{N}^K$ is $P_1 \otimes \cdots \otimes P_K$-*multi-chaotic* if, for any $m \geq 1$, the convergence in distribution

$$\lim_{N \to \infty} (X_{n,k}^N, 1 \leq n \leq m, 1 \leq k \leq K) \overset{\text{dist}}{=} P_1^{\otimes m} \otimes \cdots \otimes P_K^{\otimes m}$$

holds for the topology of uniform convergence on compact sets, where $P_k$ for $1 \leq k \leq K$ is a probability distribution on $\mathbb{R}_+$ and with the convention that $N$ goes to infinity when $\min_k N_k$ goes to infinity.

If a sequence of random variables is multi-chaotic, then any fixed finite sub-system is asymptotically independent with particles of class $k$ having the law $P_k$.

An interesting result, not used in this paper, is that a sequence of multi-exchangeable multi-class systems is multi-chaotic *if and only if* the restriction to each class is chaotic. See Theorem 3 in Graham [14].

*Mean field scaling: Initial conditions and definitions.*   Equation (3.4) defines the growth rate of the number of class $k$ connections when the total number of connections goes to infinity, one denotes

$$M^N = \min_{1 \leq k \leq K} N_k,$$

and the terminology "$N$ goes to infinity" refers to the fact that $M^N$ converges to infinity. Additionally, it will be assumed that the initial conditions of SDE (3.5),

$$(W_{n,k}^N(0), 1 \leq n \leq N_k, 1 \leq k \leq K),$$

are multi-exchangeable and $P_{1,0} \otimes \cdots \otimes P_{K,0}$-multi-chaotic where, for $1 \leq k \leq K$, $P_{k,0}$ is a probability distribution on $\mathbb{R}_+$. The process $(W_{n,k}^N(t))$ is the solution of (3.5) in this setting. For $N \in \mathbb{N}^K$ and $n \leq M^N$, one defines

$$W_n^N(t) = (W_{n,k}^N(t), 1 \leq k \leq K).$$

Note that, because of exchangeability, $(W_n^N(t), t \geq 0)$ and $(W_1^N(t), t \geq 0)$ have the same distribution, and, because of multi-chaoticity, the variable $W_1^N(0)$ converges in distribution to $P_{1,0} \otimes \cdots \otimes P_{K,0}$ as $M^N$ goes to infinity.

The following theorem is the main result for the mean-field convergence. The underlying topology on the corresponding functional space is the uniform convergence on compact sets.

**Theorem 5.1.** *Under the conditions:*

(1) *for $1 \leq k \leq K$,*

$$\lim_{N \to +\infty} \frac{N_k}{N_1 + \cdots + N_K} = p_k;$$

(2) *the variables $(W_{n,k}^N(0), 1 \leq n \leq N_k, 1 \leq k \leq K)$ are multi-exchangeable and $P_{1,0} \otimes \cdots \otimes P_{K,0}$-multi-chaotic;*

(3) *the functions $a_k \colon \mathbb{R}_+ \times \mathbb{R}_+^J \to \mathbb{R}_+$, $1 \leq k \leq J$, are bounded and Lipschitz, and Condition (C) holds for the random variables $\{W_1^N(0), N \in \mathbb{N}^K\}$, $(b_k)$ and $\varepsilon > 0$;*

*then, for $1 \leq k \leq K$,*

$$\lim_{N \to +\infty} \mathbb{E}\left| \frac{1}{N_k} \sum_{n=1}^{N_k} W_{n,k}^N(t) - \mathbb{E}(W_k(t)) \right| = 0,$$

*and the sequence of processes*

$$((W_{n,k}^N(t), t \geq 0), 1 \leq n \leq N_k, 1 \leq k \leq K),$$

*and the solutions of SDE (3.5) with initial conditions $(W_{n,k}^N(0))$, is multi-exchangeable and $P_W$-multi-chaotic where $P_W = P_{W_1} \otimes \cdots \otimes P_{W_K}$ is the distribution of the process $(W(t), t \geq 0) = ((W_k(t), t \geq 0), 1 \leq k \leq K)$, the solution of the nonlinear SDE (4.1) with initial distribution $P_{1,0} \otimes \cdots \otimes P_{K,0}$.*

*In particular, for any $n \geq 1$, $(W_n^N(t), t \geq 0) = ((W_{n,k}^N(t), 1 \leq k \leq K), t \geq 0)$ converges in distribution to $(W(t), t \geq 0)$ as $N$ goes to infinity.*

**Proof.** We prove the result assuming that $b_k$ is Lipschitz. The other case is similar (see, e.g., the end of the proof of Theorem 4.2). A coupling method will be used.

Let $(V_n^N(t)) = ((V_{n,k}^N(t), 1 \leq k \leq K))$ be the solution of (4.1) with the initial conditions of the process $(W_n^N(t))$,

$$V_n^N(0) = (V_{n,k}^N(0), 1 \leq k \leq K) = (W_{n,k}^N(0), 1 \leq k \leq K) = W_n^N(0)$$

and Poisson processes $(\mathcal{N}_{n,k}, 1 \leq k \leq K)$ driving the dynamic are the same as for $(W_n^N(t))$ in (3.5).

The key result to prove the theorem is to prove that the processes $(W_n^N(t))$ and $(V_n^N(t))$ are asymptotically close, that is, for any $T > 0$ and $n \geq 1$,

$$(5.1) \qquad \lim_{N \to +\infty} \mathbb{E}\left( \sup_{0 \leq t \leq T} \|W_n^N(t) - V_n^N(t)\| \right) = 0.$$

Because of the multi-exchangeability property, one has to consider only the case $n = 1$. With the notation of Proposition 4.1, one has $V_1^N = \phi(W_1^N(0), u_{V_1^N})$ by definition of $V_1^N$, and $W_1^N = \phi(W_1^N(0), \overline{U}^N(t))$ by (3.5). By taking $X_0 = X_0' = W_1^N(0)$, $u = u_{V_1^N}$, and $u' = \overline{U}^N$, and the expectation of inequality (4.3) with respect to the Poisson processes $(\mathcal{N}_{n,k}, 1 \leq k \leq K, n \neq 1)$, gives

$$
\begin{aligned}
(5.2) \qquad &\mathbb{E}^N\Big(\sup_{s \leq t} \|W_1^N(s) - V_1^N(s)\|\Big) \\
&\leq C_1^N(t) e^{C_1^N(t)t} \int_0^t \mathbb{E}^N(\|\overline{U}^N(s) - u_{V_1^N}(s)\|)\,ds,
\end{aligned}
$$

where $\mathbb{E}^N(\cdot)$ denotes the expectation conditionally on the random vector $W_1^N(0)$, and, since the function $(a_k)$ is bounded by some $\|a\|$,

$$(5.3) \quad V_{1,k}^N(t) \leq V_{1,k}^N(0) + t\|a\| = W_{1,k}^N(0) + t\|a\|, \qquad 1 \leq k \leq K, t \geq 0,$$

and $C_1^N(t)$ can be taken as $A(1 + \|W_1^N(0)\| + t)$, for some constant $A$ depending only on the functions $(a_k)$ and $(b_k)$.

We hereafter cope with the lack of symmetry of multi-exchangeable systems, whereas similar proofs exploit the strength of full exchangeability.

It is assumed that $p_k > 0$ for all $1 \leq k \leq K$ for ease of presentation (the proof extends easily), then, for $t \geq 0$,

$$(5.4) \quad \mathbb{E}^N(\|\overline{U}^N(t) - u_{V_1^N}(t)\|) \leq \sum_{j=1}^J \sum_{k=1}^K A_{jk} \mathbb{E}^N\left|\frac{N_k}{|N|}\overline{W}_k^N(t) - p_k\mathbb{E}(V_{1,k}^N(t))\right|$$

with, for $1 \leq k \leq K$,

$$\overline{W}_k^N(t) = \frac{1}{N_k}\sum_{n=1}^{N_k} W_{n,k}^N(t)$$

and

$$
\begin{aligned}
(5.5) \qquad &\mathbb{E}^N\left|\frac{N_k}{|N|}\overline{W}_k^N(t) - p_k\mathbb{E}(V_{1,k}^N(t))\right| \\
&\leq p_k\mathbb{E}^N\left|\frac{1}{N_k}\sum_{n=1}^{N_k}[W_{n,k}^N(t) - \mathbb{E}(V_{1,k}^N(t))]\right| + \Delta_k^N(t)
\end{aligned}
$$

with

$$\Delta_k^N(t) \stackrel{\text{def}}{=} \left|\frac{1}{|N|} - \frac{p_k}{N_k}\right| \mathbb{E}^N\left|\sum_{n=1}^{N_k} W_{n,k}^N(t)\right| \leq \left|\frac{N_k}{|N|} - p_k\right| \mathbb{E}^N|W_{1,k}^N(t)|,$$

and, consequently,

$$\lim_{N \to \infty} \mathbb{E}\Big(\sup_{s \leq t} \Delta_k^N(s)\Big) = 0$$

due to the uniform integrability bounds on the initial data and the growth bound of the solutions of (3.5).

The finite set $\{1, \ldots, N_k\}$ is partitioned into $\lfloor N_k/M^N \rfloor$ consecutive subsets with $M^N$ elements and the remaining subset with cardinality $R_k^N < M^N$. Using again the multi-exchangeability property, one gets

$$
\begin{aligned}
\mathbb{E}^N &\left| \frac{1}{N_k} \sum_{n=1}^{N_k} [W_{n,k}^N(t) - \mathbb{E}(V_{1,k}^N(t))] \right| \\
(5.6) \qquad &\leq \left\lfloor \frac{N_k}{M^N} \right\rfloor \frac{M^N}{N_k} \mathbb{E}^N \left| \frac{1}{M^N} \sum_{n=1}^{M^N} [W_{n,k}^N(t) - \mathbb{E}(V_{1,k}^N(t))] \right| \\
&\quad + \frac{R_k^N}{N_k} \mathbb{E}^N \left| \frac{1}{R_k^N} \sum_{n=1}^{R_k^N} [W_{n,k}^N(t) - \mathbb{E}(V_{1,k}^N(t))] \right|,
\end{aligned}
$$

where

$$
\begin{aligned}
\mathbb{E}^N &\left| \frac{1}{M^N} \sum_{n=1}^{M^N} [W_{n,k}^N(t) - \mathbb{E}(V_{1,k}^N(t))] \right| \\
(5.7) \qquad &\leq \mathbb{E}^N \left| \frac{1}{M^N} \sum_{n=1}^{M^N} (W_{n,k}^N(t) - V_{n,k}^N(t)) \right| \\
&\quad + \mathbb{E}^N \left| \frac{1}{M^N} \sum_{n=1}^{M^N} (V_{n,k}^N(t) - \mathbb{E}(V_{1,k}^N(t))) \right| \\
&\leq \mathbb{E}^N |W_{1,k}^N(t) - V_{1,k}^N(t)| + \mathbb{E}^N \left| \frac{1}{M^N} \sum_{n=1}^{M^N} [V_{n,k}^N(t) - \mathbb{E}(V_{1,k}^N(t))] \right|.
\end{aligned}
$$

The Jensen inequality yields

$$
\begin{aligned}
\mathbb{E} &\left( \left( \mathbb{E}^N \left| \frac{1}{M^N} \sum_{n=1}^{M^N} [V_{n,k}^N(t) - \mathbb{E}(V_{1,k}^N(t))] \right| \right)^2 \right) \\
&\leq \mathbb{E} \left( \left| \frac{1}{M^N} \sum_{n=1}^{M^N} [V_{n,k}^N(t) - \mathbb{E}(V_{1,k}^N(t))] \right|^2 \right),
\end{aligned}
$$

and by developing the square, one gets

$$
\begin{aligned}
\mathbb{E} &\left( \left| \frac{1}{M^N} \sum_{n=1}^{M^N} [V_{n,k}^N(t) - \mathbb{E}(V_{1,k}^N(t))] \right|^2 \right) \\
&= \frac{1}{M^N} \mathbb{E}([V_{1,k}^N(t) - \mathbb{E}(V_{1,k}^N(t))]^2)
\end{aligned}
$$

$$+ \frac{M^N - 1}{M^N} \mathbb{E}([V_{1,k}^N(t) - \mathbb{E}(V_{1,k}^N(t))][V_{2,k}^N(t) - \mathbb{E}(V_{1,k}^N(t))]),$$

and by exchangeability, the covariance term of the last expression is given by

$$\mathbb{E}(V_{1,k}^N(t)V_{2,k}^N(t)) - 2\mathbb{E}(V_{1,k}^N(t))\mathbb{E}(V_{1,k}^N(t)) + [\mathbb{E}(V_{1,k}^N(t))]^2.$$

Hence the uniform integrability bounds yield

$$\lim_{N \to \infty} \mathbb{E}\left(\left(\sup_{s \le t} \mathbb{E}^N \left| \frac{1}{M^N} \sum_{n=1}^{M^N} [V_{n,k}^N(s) - \mathbb{E}(V_{1,k}^N(s))] \right|\right)^2\right) = 0.$$

For the remainder term containing $R_k^N$ in (5.6), two cases have to be considered,

$$R_k^N \le \sqrt{|N|} \quad \Rightarrow \quad \lim \frac{R_k^N}{N_k} = 0 \quad \text{and} \quad R_k^N > \sqrt{|N|} \quad \Rightarrow \quad \lim R_k^N = \infty,$$

and the second case is treated as for $M^N$. Inequalities (5.4), (5.5), (5.6) and (5.7) yield

$$\mathbb{E}^N \|u_{V_1^N}(t) - \overline{U}^N(t)\| \le \sum_{j=1}^{J} \sum_{k=1}^{K} A_{jk} p_k \mathbb{E}^N \|V_1^N(t) - W_1^N(t)\| + F^N(t)$$

with

$$(5.8) \qquad\qquad \lim_{N \to \infty} \mathbb{E}\left(\sup_{s \le t} |F^N(s)|\right) = 0.$$

By plugging this bound in inequality (5.2), one gets

$$\mathbb{E}^N(\|W_1^N - V_1^N\|_t) \le C_1^N(t)e^{C_1^N(t)t}\left(A \int_0^t \mathbb{E}^N(\|W_1^N - V_1^N\|_s)\,ds + F^N(t)\right)$$

for some constant $A > 0$. The Gronwall lemma gives the relation

$$\mathbb{E}^N(\|W_1^N - V_1^N\|_t) \le F^N(t)C_1^N(t)e^{C_1^N(t)t}\exp(AtC_1^N(t)e^{C_1^N(t)t}),$$

the exponential moment assumption for the sequence $(W_1^N(0), \in \mathbb{N}^K)$ shows that, for $t < t_0$ with $t_0$ depending only on the functions $(a_k)$ and $(b_k)$ and $\varepsilon$, the random variables $C_1^N(t)\exp(C_1^N(t)t)$, $N \ge 1$, are tight. Relation (5.8) implies, therefore, that the random variable $\mathbb{E}^N(\|V_1^N - W_1^N\|_t)$ converges to 0 in probability when $N$ goes to infinity. The growth bounds (5.3) (also valid for $W_1^N$) and the uniform integrability property of $(W_1^N(0), N \in \mathbb{N}^K)$ show that the convergence to 0 also holds for the expected value by Lebesgue's theorem. Relation (5.1) has therefore been established for $t \le t_0$. The extension to any arbitrary $t$ is done as in the proof of Theorem 4.2.

The continuity property with respect to the initial conditions implied by relation (4.7) and the uniqueness result of Theorem 4.2 show that the sequence of processes $(V_{n,k}^N(t), 1 \leq n \leq N_k, 1 \leq k \leq K)$ is $P_W = P_{W_1} \otimes \cdots \otimes P_{W_K}$-multi-chaotic. Relation (5.1) gives that the same property also holds for the sequence of processes $(W_{n,k}^N(t), 1 \leq n \leq N_k, 1 \leq k \leq K)$. The theorem is proved. $\square$

The following corollary is an immediate consequence of the above result and of the exchangeability property.

COROLLARY 5.2. *Under the assumptions of Theorem 5.1 and if, for $1 \leq k \leq K$, $P_{W_k}$ is the distribution of $(W_k(t))$ such that the process $(W(t)) = ((W_k(t)), 1 \leq k \leq K)$ is the solution of the nonlinear SDE (4.1) with initial distribution $P_{1,0} \otimes \cdots \otimes P_{K,0}$, then the convergence in distribution of the empirical distributions*

$$\lim_{N \to \infty} \Lambda_k^N = P_{W_k}$$

*holds for the weak topology on $\mathcal{P}(\mathcal{D}(\mathbb{R}_+, \mathbb{R}_+))$ with $\mathcal{D}(\mathbb{R}_+, \mathbb{R}_+)$ endowed with the Skorohod topology.*

## 6. Fixed points and invariant distributions.

This section investigates the stationary properties of the solutions of stochastic differential equations (4.1). Because of Poisson processes driving the corresponding random events, a solution $(W(t))$ of the system of (4.1) with some deterministic initial state clearly has the Markov property. Since the function $(u_W(t))$ has a role in these equations, it will imply that, in general, the Markov property is *nonhomogeneous* with time. In the martingale problem formulation, the analogue of the infinitesimal generator *depends*, a priori, not only on time but also on $\mathbb{E}(W(t))$ which complicates a lot the analysis of the asymptotic behavior of the process $(W(t))$. See the remark above in Section 4. As the following theorem shows, this has important implications, in particular, for the existence and uniqueness of an invariant distribution for $(W(t))$.

THEOREM 6.1. *In the case where, for $1 \leq k \leq K$,*

$$a_k(w, u) = a_k(u) \quad and \quad b_k(w, u) = w\beta_k(u), \qquad w \geq 0, u \in \mathbb{R}_+^J,$$

*$\beta_k$ is a positive function and $a_k$ and $b_k$ are Lipschitz functions and $a_k$ is bounded, then the invariant distributions for solutions $(W(t))$ of (4.1) are in one-to-one correspondence with the solutions $u \in \mathbb{R}_+^J$ of the fixed point equation*

$$(6.1) \qquad u_j = \sum_{k=1}^K A_{jk} p_k \psi(r_k) \sqrt{\frac{a_k(u)}{\beta_k(u)}}, \qquad 1 \leq j \leq J,$$

*where*

$$\psi(r) = \sqrt{\frac{2}{\pi}} \prod_{n=1}^{+\infty} \frac{1 - r^{2n}}{1 - r^{2n-1}}.$$

*If $u^*$ is such a solution, the corresponding invariant distribution has the density $w \to \prod_{k=1}^{K} H_{r_k, \rho_k}(w_k)$ on $\mathbb{R}_+^K$ where $\rho_k = a_k(u^*)/\beta_k(u^*)$ and $H_{r,\rho}$ is defined in Proposition 2.1.*

PROOF. Let $u^* \in \mathbb{R}_+^J$ be a solution of (6.1), and let $W_0 = (W_{0,k}) \in \mathbb{R}_+^K$ be a random variable such that, for $1 \leq k \leq K$, $W_{0,k}$ has the density (2.5) with $a = a_k(u^*)$, $b = \beta_k(u^*)$ and $r = r_k$. By the invariance of this density (see Proposition 2.1), the solution $(\widetilde{W}_k(t))$ of the stochastic differential equation

$$d\widetilde{W}_k(t) = a_k(u^*)\,dt - (1 - r_k)\widetilde{W}_k(t-)\mathcal{N}_k([0, \widetilde{W}_k(t-)\beta_k(u^*)], dt)$$

with initial condition $\widetilde{W}_k(0) = W_{0,k}$, is a stationary process, in particular, for any $t \geq 0$, $\mathbb{E}(\widetilde{W}_k(t)) = \mathbb{E}(W_{0,k})$ so that, for $1 \leq j \leq J$,

$$u_{\widetilde{W},j}(t) = \sum_{k=1}^{K} A_{jk} p_k \mathbb{E}(\widetilde{W}_k(t)) = \sum_{k=1}^{K} A_{jk} p_k \psi(r_k) \sqrt{\frac{a_k(u^*)}{\beta_k(u^*)}} = u_j^*$$

by (2.6) and (6.1) so that $u_{\widetilde{W}}(t) \equiv u^*$. This implies that the stochastic process $(\widetilde{W}_k(t), 1 \leq k \leq K)$ satisfies (4.1). The variables $W_{0,k}$ having a Gaussian moment, the above theorem shows that there exists a unique solution $(W(t))$ of (4.1) which is thus stationary and such that $u_W(t) \equiv u^*$.

Conversely, if $\pi$ is a stationary distribution of $(W(t))$, then if the distribution of $W(0)$ is $\pi$, the quantity $u^* = u_W(t)$ does not depend on $t$ and, consequently, for $1 \leq k \leq K$, the $k$th marginal of $\pi$ is the invariant distribution of the Markov process whose infinitesimal generator is given by

$$\Omega(f)(x) = a_k(u^*)f'(x) + x\beta_k(u^*)(f(rx) - f(x)).$$

By Proposition 2.1, the density of the $k$th marginal of $\pi$ is necessarily $H_{r_k, \rho_k}$ with $\rho_k = a_k(u^*)/\beta_k(u^*)$, and $u^*$ satisfies (6.1). The proposition is proved. □

The above theorem shows that if the fixed point equation (6.1) has several solutions, then the limiting process $(W(t))$ has several invariant distributions. Similarly, if (6.1) has no solution, then, in particular, $(W(t))$ cannot converge to an equilibrium. These possibilities have been suggested in the Internet literature through simulations, like the cyclic behavior of some nodes in the case of congestion. Obviously, there exist coefficients $(a_k)$ and $(\beta_k)$ with these possibilities for (6.1). Under mild and natural assumptions, such that if the loss rate at any node is increasing with its utilization, for some networks there exists a unique fixed point. See the examples investigated below.

Several special cases are now considered under the assumptions of Theorem 6.1. It is assumed that $A_{jk} = 1$ if node $j \in \{1, \ldots, J\}$ is used by connections of class $k$ and 0 otherwise. In the same way, for $k \in \{1, \ldots, K\}$, the function $\beta_k$ for the loss rate for class $k$ connections depends only on the sum of the utilizations of the nodes used by these connections, that is,

$$\beta_k(u) = \beta_k\left(\sum_{j=1}^{J} A_{jk} u_j\right)$$

with a slight abuse of notation. Equation (6.1) becomes, in this context,

$$u_j = \sum_{k=1}^{K} A_{jk} p_k \psi(r_k) \frac{\sqrt{a_k}}{\sqrt{\beta_k(\sum_{j=1}^{J} A_{jk} u_j)}}, \qquad 1 \le j \le J.$$

It is assumed that $u \to \beta_k(u)$ is strictly nondecreasing and Lipschitz and that the function $(a_k(u))$ is constant equal to $a_k$. Let $\rho_k(u)$ denote $a_k/\beta_k(u)$.

*One single node.*  Here $J = 1$, for $k \in \{1, \ldots, K\}$, there is a unique solution $u = U$ to (6.1),

$$u = \sum_{k=1}^{K} \psi(r_k) p_k \sqrt{\frac{a_k}{\beta_k(u)}}$$

and consequently a unique invariant distribution for which the state of a class $k$ connection has the density $H_{r,\rho_k}$ and its expected value, the throughput, is

$$\psi(r_k) \sqrt{\frac{a_k}{\beta_k(U)}}.$$

In particular, if the effects of the congestion of the node on a connection are independent of its class, that is, $\beta_k \equiv \beta$ for $1 \le k \le K$, then, if the scaling factors $(r_k)$ are equal, the throughputs of the different classes of connections differ by the multiplicative factors $(\sqrt{a_k})$. Recall that $a_k$ is typically proportional to the inverse of the round-trip time of a class $k$ connection. See Figure 1. See Section 3.

*A linear network.*  A network with $J$ nodes and $J + 1 = K$ classes of connections is considered. For $1 \le j \le J$, connections of class $j$ uses only node $j$ and the route of connections of class $J + 1$ go along the $J$ nodes. This network has been considered in Ben Fred et al. [5]. In this context, (6.1) becomes

$$(6.2) \qquad u_j = \left(\frac{\alpha_j}{\sqrt{\beta_j(u_j)}} + \frac{\alpha_{J+1}}{\sqrt{\beta_{J+1}(\|u\|)}}\right), \qquad 1 \le j \le J,$$

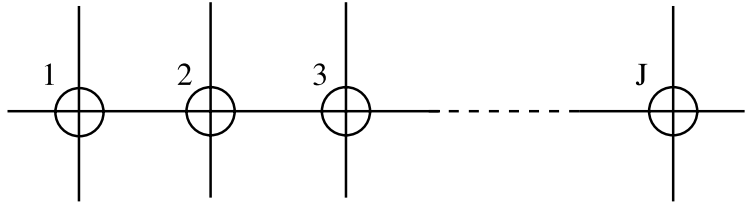

Fig. 1. *A network with $J$ nodes and $J+1 = K$ classes of connections.*

with $\alpha_j = \psi(r_j) p_j \sqrt{a_j}$, $1 \le j \le J+1$ and $\|u\| = |u_1| + \cdots + |u_J|$. For $\alpha > 0$, the function $\phi_{j,\alpha}(x) \stackrel{\text{def}}{=} x - \alpha/\sqrt{\beta_j(x)}$ being strictly nondecreasing, the equation

$$u = \sum_{j=1}^{J} \phi_{j,\alpha_j}^{-1}\left(\frac{\alpha_{J+1}}{\sqrt{\beta_{J+1}(u)}}\right)$$

has a unique solution $U$. A solution $(u_j)$ of (6.2) is necessarily given by

$$u_j = \phi_{j,\alpha_j}^{-1}\left(\frac{\alpha_{J+1}}{\sqrt{\beta_{J+1}(U)}}\right), \qquad 1 \le j \le J.$$

Equation (6.2) has, therefore, a unique solution.

*Network on a torus.* For this network, there are $J$ nodes and $K = J$ classes of connections, and class $j \in \{1, \ldots, J\}$ uses two nodes: node $j$ and $j+1$. The assumptions on the rate functions are the same as before. See Figure 2. For this network, (6.1) becomes

$$(6.3) \quad u_j = \left(\frac{\alpha_{j-1}}{\sqrt{\beta_j(u_{j-1} + u_j)}} + \frac{\alpha_j}{\sqrt{\beta_{j+1}(u_j + u_{j+1})}}\right), \qquad 1 \le j \le J,$$

with the convention that $J+1 = 1$ and $0 = J$ for the subscripts. For $1 \le j \le J$, define $\alpha_j = \psi(r_j) p_j \sqrt{a_j}$ and $\phi_j(x) = 1/\sqrt{\beta_j(x)}$. If $y_j = \alpha_{j-1}\phi_j(u_{j-1} + u_j)$, then the above equation can be rewritten as follows:

$$(6.4) \quad y_j = \alpha_{j-1}\phi_j(y_{j-1} + 2y_j + y_{j+1}), \qquad 1 \le j \le J.$$

Take $J = 3$. Since, for $1 \le j \le J$, $\phi_j$ is strictly decreasing and continuous, for $s > 0$ and $\alpha > 0$, there exists unique $x_{j,\alpha}(s)$ such that, $x_{j,\alpha}(s) = \alpha \phi_j(x_{j,\alpha}(s) + s)$ and the function $s \to x_{j,\alpha}(s)$ is strictly decreasing and continuous on $\mathbb{R}_+$. Similarly there exists a unique $S > 0$ satisfying the relation

$$S = x_{1,\alpha_1}(S) + x_{2,\alpha_2}(S) + x_{3,\alpha_3}(S).$$

From this one sees that $y = (x_{i,\alpha_i}(S), 1 \le i \le 3)$ is the unique solution of (6.4). One concludes that (6.3) has also a unique solution.

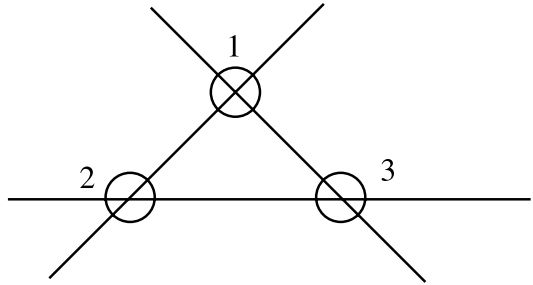

Fig. 2. *Triangular network.*

UMR 7641 CNRS—ÉCOLE POLYTECHNIQUE
ROUTE DE SACLAY
91128 PALAISEAU
FRANCE
E-MAIL: carl@cmapx.polytechnique.fr

INRIA PARIS—ROCQUENCOURT
DOMAINE DE VOLUCEAU
78153 LE CHESNAY
FRANCE
E-MAIL: Philippe.Robert@inria.fr
URL: http://www-rocq.inria.fr/˜robert